\documentclass[11pt]{article}

\usepackage[dvips]{graphicx}
\usepackage{amssymb}

\usepackage{caption}
\usepackage{hyperref}
\usepackage{arcs}
\usepackage {xcolor} 
\usepackage{subfigure}

\newtheorem{theorem}{Theorem}
\newtheorem{Proposition}{Proposition}

\newtheorem{Remark}{Remark}

\newenvironment{AMS}{\small\bf 2010 AMS subject classification: }{} 

\begin{document}

\title{Qsurf: compressed QMC integration\\ on parametric surfaces}
\author{G. Elefante, A. Sommariva, M. Vianello \\ \\
University of Padova}

\date{\today}  

 \maketitle

\begin{abstract} 
We discuss a ``bottom-up'' algorithm for Tchakaloff-like compression of Quasi-MonteCarlo (QMC) integration on surfaces that admit an analytic parametrization. The key tools are  
Davis-Wilhelmsen theorem on the so-called “Tchakaloff sets” for positive linear functionals on polynomial spaces, and Lawson-Hanson algorithm for NNLS. This algorithm shows remarkable speed-ups with respect to Caratheodory-like subsampling, since it is able to work with much smaller matrices. We provide the corresponding Matlab code {\em Qsurf}, together with integration tests on regions of different surfaces such as sphere and torus. 
\end{abstract}

\vskip0.2cm 
\noindent
\begin{AMS}
{\rm 65C05, 65D32.}
\end{AMS}

\section{Introduction}

In the recent paper \cite{ESV23}, we have considered the compression problem for Quasi-MonteCarlo (QMC) surface integration on multibubbles (the surface of a ball union in $\mathbb{R}^3$), which can have a quite complicated structure. Indeed, numerical modelling with multibubbles is relevant in several applications, but compression of QMC integration seemed an overlooked approach, especially in the case of surface integrals. 

In this paper, we extend such an approach to compressed QMC formulas for  general integration problems on compact subsets of surfaces in $\mathbb{R}^3$, admitting an analytic parametrization. Such formulas preserve the approximation power of QMC up to the best uniform polynomial approximation error of a given degree to the integrand, but using a much lower number of sampling points. 

The key tools are Davis-Wilhelmsen theorem on the so-called ``Tchakaloff sets'' for positive linear functionals and Lawson-Hanson algorithm for 
NNLS, which allows to extract a set of ``equivalent'' re-weighted nodes from a huge uniformly distributed sequence with respect to the surface measure, by working in a ``bottom-up'' mode. Such a sequence can be obtained for example from a bivariate Halton sequence by an area-preserving map, when available, or by the probabilistic method of rejection sampling, which has been extended to the low-discrepancy deterministic setting, cf. e.g. \cite{NO16,ZD14}. On the other hand, there are 
other relevant QMC point sequences on manifolds, see e.g. \cite{BCCGST14,BSSW14}.

The ``bottom-up'' approach  shows remarkable speed-ups with respect to Caratheodory-like subsampling (cf. e.g. \cite{H21,LL12,PSV17,SV15,Tche15}), 
since it is able to work with much smaller matrices. We stress that one of the main difficulties consists in adapting the compression algorithm
to work on the appropriate spaces of trivariate polynomials restricted to the surface, since the dimension of trivariate polynomial spaces can collapse 
in the case of algebraic surfaces.

The paper is organized as follows. In Section 2 we briefly discuss the theoretical
background and the main idea of the ``bottom-up'' compression algorithm. Then, we sketch the algorithm, that has been implemented in Matlab, and comment on the main computational issues. Finally, in Section 3 we present some numerical examples concerning regions of sphere and torus, and the Cartesian graph of an analytic function. All the codes and demos are all freely available at \cite{ESV23-2}.

\section{QMC compression on surfaces}
The possibility of compressing QMC integration rests on a somehow overlooked but relevant result of quadrature theory, originally proved by Davis \cite{D67} and then extended by Wilhelmsen \cite{W76}. Only recently this theorem has been rediscovered as a basic tool for positive cubature via adaptive NNLS moment-matching, cf. \cite{ESV22,L21,SV21,SV23}.

\begin{theorem} (Davis, 1967 - Wilhelmsen, 1976) Let $\{f_j\}_{1\leq j\leq N}$ be continuous, 
real-valued, linearly independent functions 
defined on a compact set $\Omega \subset {\mathbb{R}}^d$, 
and $\mathcal{F}=span(f_1,\dots,f_N)$. 
Assume that $\mathcal{F}$ satisfies the Krein condition (i.e. there is at least 
one $f \in \mathcal{F}$ which does not vanish on $\Omega$) and that $L$ is a positive 
linear functional on $\mathcal{F}$, i.e. $L(f)>0$ for every $f\in \mathcal{F}$, $f \geq 0$ not vanishing everywhere in $\Omega$.  

If $\{P_i\}_{i=1}^{\infty}$ is an everywhere dense subset of $\Omega$, 
then for sufficiently large $m$, the set 
$X_m=\{P_i\}_{i=1,\ldots,m}$ is 
a ``Tchakaloff set'', i.e. there exist weights $w_k > 0$, $k=1,\ldots,\nu$, and nodes $\{Z_k\}_{k=1,\dots,\nu}
\subset X_m \subset \Omega$, 
with $\nu={\mbox{card}}(\{Z_k\}) \leq N$, such that 
\begin{equation} \label{tchq}
L(f)=\ell(f)=\sum_{k=1}^\nu w_k f(Z_k)\;,\;\;\forall f \in \mathcal{F}\;.
\end{equation}
\end{theorem}

\vskip0.5cm

Davis-Wilhelmsen theorem is a constructive generalization of the well-known Tchakaloff theorem \cite{T57} on the existence of positive quadrature formulas. But, just in view of its generality, it can 
be directly applied to a discrete functional like a QMC formula on $\Omega= \mathcal{J}$, $\mathcal{J}$ being a compact region  of a surface $\mathcal{S}\subset \mathbb{R}^3$
\begin{equation} \label{QMC}
L(f)=L_{\mbox{\tiny{QMC}}}(f)=\frac{\sigma(\mathcal{J})}{M}\,\sum_{i=1}^M{f(P_i)}\approx 
\int_\mathcal{J}{f}\,d\sigma\;,\;\;f\in C(\mathcal{J})\;,
\end{equation}
where 
$$X_M=\{P_i\}_{i=1,\ldots,M}\;,\;\;M>N\;,$$ 
is a uniformly distributed sequence on $\mathcal{J}$ 
and $\sigma$ is the surface measure. 
Typically one generates a uniformly distributed sequence of cardinality say $M_0$ on the bounding surface $\mathcal{S}\supseteq \mathcal{J}$, from which sequence on $\mathcal{J}$ is extracted by a suitable in-domain algorithm.  We observe that if $\sigma(\mathcal{J})$ is unknown or difficult to compute, it can be approximated  as $\sigma(\mathcal{J})\approx \sigma(\mathcal{S})M/M_0$.

Positivity of the functional for $f\in \mathcal{F}=\mathbb{P}_n^3(\mathcal{J})$ 
(the space of trivariate polynomials of total degree not exceeding $n$ restricted to $\mathcal{J}$), is ensured whenever the set $X_M$ is $\mathbb{P}_n^3(\mathcal{J})$-determining, i.e. a polynomial vanishing there vanishes everywhere on $\mathcal{J}$, 
or equivalently $dim(\mathbb{P}_n^3(X_M))=N=dim(\mathbb{P}_n^3(\mathcal{J}))$, or even 
\begin{equation} \label{ranks}
rank(V_M)=N\;,\;\;V_M=V^{(n)}(X_M)=[f_j(P_i)]\in \mathbb{R}^{M\times N}\;
\end{equation}
where $V_M$ is the corresponding rectangular Vandermonde-like matrix. Notice that, $X_M$ being a sequence, for every $k\leq M$ we have that 
\begin{equation} \label{submatrices}
V_k=V^{(n)}(X_k)=[(V_M)_{ij}]\;,\;\;1\leq i \leq k\;,\;1\leq j\leq N\;.
\end{equation}

We stress that the full rank requirement for $V_M$ is not restrictive, in practice, when $\mathcal{S}$ is a surface that admits an analytic parametrization, the subset $\mathcal{J}$ is $\mathbb{P}_n^3(\mathcal{S})$-determining and the points are uniformly distributed with respect to the surface measure. 
Indeed, the probability that $det(V_N)=0$ dealing with uniformly distributed points is null, as is ensured by the following proposition which is a special case of a general result proved in \cite{DASV23} in the case of continuous random point distributions. 

\begin{Proposition}
Let $\mathcal{S}\subset \mathbb{R}^3$ be a surface that admits an analytic parametrization $P=\Psi(u,v)$ from a connected open set $D\subset \mathbb{R}^2$, 
i.e. $\Psi=(\Psi_1,\Psi_2,\Psi_3)$ where $\Psi_i:D\to \mathbb{R}^3$ are analytic 
and $\Psi(D)=\mathcal{S}$. Moreover, let $\{f_j\}_{1\leq j\leq N}$ be a basis of $\mathbb{P}_n^3(\mathcal{S})$ and $\{(u_i,v_i)\}_{i\geq 1}$ an  equidistributed sequence on $D$ with respect to any given probability density $\phi(u,v)$. 

Then, the points $\{P_i=\Psi(u_i,v_i)\}_{1\leq i\leq N}$ are almost surely unisolvent 
for polynomial interpolation in $\mathbb{P}_n^3(\mathcal{S})$.
\end{Proposition}
\vskip0.2cm

\begin{Remark}
{\em 
We can apply this proposition to the case where the parametrization is regular (so that the surface area element $\|\partial_u\Psi\times \partial_v\Psi\|_2/\sigma(\mathcal{J})$ is well-defined), $dim(\mathbb{P}_n^3(\mathcal{S}))=N
=dim(\mathbb{P}_n^3(\mathcal{J}))$, and $d\sigma=\phi(u,v)\,du\,dv$ with density
\begin{equation} \label{surfdens}
\phi(u,v)=I_{{\mathcal{J}}}(\Psi(u,v))\,\|\partial_u\Psi\times \partial_v\Psi\|_2/\sigma(\mathcal{J})\;,
\end{equation}
$I_{\mathcal{J}}$ denoting the indicator function of $\mathcal{J}$.
}
\end{Remark}

\begin{Remark} 
{\em 
To be rigorous, we should notice that Proposition 1 concerns random sequences, whereas here we deal with quasi-random sequences, where we can expect, and we have indeed verified experimentally, that the full-rank property of $V_M$ in practice holds. In order to construct a sequence $\{P_i\}_{1\leq i\leq N}$ that be uniformly distributed on $\mathcal{J}$ with respect to the surface measure, we can adopt the classical 
probabilistic method of {\em rejection sampling} on $D$ applied to the density (\ref{surfdens}), that has been extended to low-discrepancy sequences; cf. \cite{NO16,ZD14} with 
the references therein. Clearly, a suitable ``in-domain'' algorithm for $\mathcal{J}$ has to be at hand. 
}
\end{Remark}

\begin{Remark} 
{\em 
We recall that polynomial spaces can collapse on algebraic surfaces, i.e. it happens that $dim(\mathbb{P}_n^3(\mathcal{J}))=dim(\mathbb{P}_n^3(\mathcal{S}))<dim(\mathbb{P}_n^3(\mathbb{R}^3))=(n+1)(n+2)(n+3)/6$. For example, if $\mathcal{J}$ is a subset with internal points 
w.r.t. the topology of the sphere $S^2$ (e.g. a spherical polygon as in the first example below), we have that $dim(\mathbb{P}_n^3(\mathcal{J}))=dim(\mathbb{P}_n^3(S^2))=(n+1)^2$; we refer the reader, e.g., to \cite{CLOS15} concerning the delicate matter of determining polynomial spaces dimension on algebraic varieties. 
}
\end{Remark}

   In view of the results quoted above, when $M\gg N$ we can then try to find a Tchakaloff set $X_m\subset X_M$, with $N\leq m<M$, such that there exists a sparse nonnegative solution vector $u$ to the underdetermined {\em moment-matching system}
\begin{equation} \label{mom-match}
V_m^t u=\lambda=V_M^t e\;,\;\;e=\frac{\sigma(\mathcal{J})}{M}\,(1,\dots,1)^t\;.
\end{equation}

In practice, we solve (\ref{mom-match}) via Lawson-Hanson active-set method \cite{LH95} applied to the NNLS problem 
\begin{equation} \label{NNLS}
\min_{u\geq 0} \|V_m^t u-\lambda\|_2\;, 
\end{equation}
accepting the solution when the residual size is small, say 
\begin{equation} \label{res-tol}
\|V_m^t u-\lambda\|_2<\varepsilon
\end{equation}
where $\varepsilon$ is a given tolerance. Then the nonzero components of $u$  provide nodes and weights of a compressed QMC formula extracted from $X_m$, that is $\{w_k\}=\{u_i:\,u_i>0\}$ and $\{Z_k\}=\{P_i:\,u_i>0\}$, giving 
\begin{equation} \label{compressed}
\ell_{\mbox{\tiny{QMC}}}(f)=\sum_{k=1}^\nu w_k f(Z_k)\;,\;\;\nu \leq N\ll M\;\;,
\end{equation} 
where $\ell_{\mbox{\tiny{QMC}}}(f)=L_{\mbox{\tiny{QMC}}}(f)$ for every $f\in \mathbb{P}_n^3(\mathcal{J})$.

It is worth recalling that, in the case $m=M$, Caratheodory theorem on finite-dimensional conic combinations (applied to the columns of $V_M^t$) would ensure directly the existence of a Tchakaloff-like representation of the QMC functional  
(cf. \cite{PSV17} for a discussion on this point in the general framework of discrete measure compression by ``Caratheodory-Tchakaloff subsampling''). In such a way, however, working with say an order of $10^5-10^6$ nodes, we would have to manage a huge matrix, that is we would have to solve the huge NNLS problem 
\begin{equation} \label{big-NNLS}
\min_{u\geq 0} \|V_M^t u-\lambda\|_2\;. 
\end{equation}

On the contrary, we can substantially reduce the computation cost by solving an increasing sequence of much smaller problems like (\ref{NNLS}) with $m:=m_1,m_2,m_3,\dots$, $m_1<m_2<m_3<\dots\leq M$,
\begin{equation} \label{small-NNLS}
\min_{u\geq 0} \|V_{m_j}^t u-\lambda\|_2\;,\;\;j=1,2,3,\dots\;,\;m_1\geq N\;,
\end{equation}
corresponding to increasingly dense subsets $X_{m_1}\subset X_{m_2} \subset \dots \subseteq  X_M$, until the residual becomes sufficiently low. We may call this procedure a ``bottom-up'' approach to QMC compression. Indeed, as shown in \cite{ESV22}, 
with a suitable choice of the sequence $\{m_j\}$ the residual becomes extremely small in few iterations, with a substantial speed-up with respect to (\ref{big-NNLS}).

Now, following \cite{ESV22} it is easy to derive the following error estimate 
$$
|\ell_{\mbox{\tiny{QMC}}}(f)-\int_{\mathcal{J}}{f}\,d\sigma |
\leq \mathcal{E}_{\mbox{\tiny{QMC}}}(f) 
+2\,\mu(\mathcal{J})\,E_n(f;X_M)
$$
\begin{equation} \label{equiv}
\leq \mathcal{E}_{\mbox{\tiny{QMC}}}(f) 
+2\,\mu(\mathcal{J})\,E_n(f;\mathcal{J})\;,
\end{equation}
valid for every $f\in C(\mathcal{J})$, where 
$\mathcal{E}_{\mbox{\tiny{QMC}}}(f)=|L_{\mbox{\tiny{QMC}}}(f)-\int_{\mathcal{J}}{f}\,d\sigma|$ and 
we define $E_n(f;K)=\inf_{p\in \mathbb{P}_n^3(K)}
{\|f-p\|_{\infty,K}}$ 
with $K$ discrete or continuous compact set.

The meaning of (\ref{equiv}) is that the compressed QMC functional $\ell_{\mbox{\tiny{QMC}}}$ retains the approximation power of the original QMC formula, up to a quantity proportional to the best polynomial approximation error to $f$ in the uniform norm on $X_M$ (and hence by inclusion in the uniform norm on $\mathcal{J}$). We recall that the latter can be estimated depending on the regularity of $f$ by multivariate Jackson-like theorems, cf. e.g. \cite{P09} for volume integrals where $\mathcal{J}$ is the closure of a bounded open set 
and \cite{R71} for the case of the sphere. On the other hand, we do not deepen here the topic of QMC convergence and error estimates, in particular on manifolds, referring the reader to specific papers and monographs, like e.g. \cite{BCCGST14,BSSW14,DP10} .

\subsection{Algorithm description and computational issues} 
In this section we sketch the method implementation in the form of a pseudo-code and discuss its main computational features. 

\vskip0.8cm
\noindent
{\bf Algorithm Qsurf}: {\em Bottom-up compression of QMC integration on a compact subset $\mathcal{J}$ of a surface $\mathcal{S}\subset \mathbb{R}^3$ 
with a regular analytic parametrization on a domain $D\subset \mathbb{R}^2$}

\begin{itemize}
    
\item {\tt input}: 
the bounding surface measure $\sigma(\mathcal{S})$, possibly the measure $\sigma(\mathcal{J})$, the cardinality $M_0$ of a uniformly distributed sequence on $\mathcal{S}$,
the cardinality increase factor $\theta>1$, the moment-matching tolerance $\varepsilon$, the residual decrease threshold $\tau>1$

\item[$(i)$] generate $M_0$ low-discrepancy points on the bounding surface $\mathcal{S}\supseteq \mathcal{J}$ (for example by rejection sampling on $D$ w.r.t. the surface measure density) and extract the 
points $X=X_M=\{P_i\}_{i=1,\dots,M}$ that lie on $\mathcal{J}$ (by a suitable 
``in-domain'' algorithm)

\item[$(ii)$] if unknown, approximate $\sigma(\mathcal{J})$ as $\sigma(\mathcal{J}):=\sigma(\mathcal{S})\,M/M_0$

\item[$(iii)$] \% {\em selecting a basis of $\mathbb{P}_n^3(X)$}
    
 $(iii1)$ take a polynomial basis $\{p_1,\dots,p_{\cal{V}}\}$ of $\mathbb{P}_n^3$, ${\cal{V}}=\frac{(n+1)(n+2)(n+3)}{6}$
 
 $(iii2)$ compute the Vandermonde-like matrix $C:=[p_j(P_i)]\in \mathbb{R}^{M\times {\cal{V}}}$
 
 $(iii3)$ 
 compute $N:=rank(C_{\mathcal{V}})$ where  $C_{\mathcal{V}}=[(C)_{ij}]$, $1\leq i,j\leq \mathcal{V}$
 
 $(iii4)$ compute the $QR$ factorization with column pivoting $C_{\mathcal{V}}^\pi=QR$ 
 where $\pi=(\pi_1,\dots,\pi_{\mathcal{V}})$ is the column permutation vector
 
 $(iii5)$ set $V_M:=[(C)_{ij}]$, $1\leq i\leq M$, $j=\pi_1,\dots,\pi_N$

\item[$(iv)$] 
compute the QMC moments $\lambda:=V_M^t e$, $e=\sigma(\mathcal{J})/M\,(1,\dots,1)^t$

\item[$(v)$] 

\% {\em bottom-up QMC compression}

$(v1)$ inizialize $m$, $N\leq m\ll M$ and $momtype:=0$

$(v2)$ set $V_m: =[(V_M)_{ij}]$, $1\leq i \leq m$, 
$1\leq j\leq N$

$(v3)$ compute the $QR$ factorization $V_m=Q_mR_m$

$(v4)$ if $momtype=0$ {\tt then}
 \begin{itemize}
\item compute the modified QMC moments $q_m=(R_m^{-1})^t\lambda$ by solving the system $R_m^tq_m=\lambda$ via Gaussian elimination with row pivoting
\item set $A_m=Q_m$
 \end{itemize}
{\hspace{0.85cm}}{\tt else}
 \begin{itemize}
 \item compute the modified QMC moments $q_m=(R_m^{-1})^t\lambda=(V_M R_m^{-1})^te$ 
 as $q_m:=A_M^te$, by solving the matrix equation $R_m^tA_M^t=V_M^t$ via Gaussian elimination with row pivoting
 \item set $A_m=[(A_M)_{i,j}]$, $1\leq i \leq m$, 
$1\leq j\leq N$
\end{itemize} 

$(v5)$ compute a sparse solution $u$ to the NNLS problem $$\min_{u\geq 0} \|A_m ^t u-q_m\|_2$$ by Lawson-Hanson active-set algorithm

$(v6)$ compute the relative residual $res:=\|V_m^t u-\lambda\|_2/\|\lambda\|_2$

$(v7)$ {\tt if} $res_0/res>\tau\;\&\;m<M$ {\tt then}\\
\phantom{SP} $(v7a)$ {\tt if} $momtype=0$ {\tt then} \\
\phantom{SP} \phantom{SP} - set $momtype:=1$ and {\tt goto} $(v2)$\\
\phantom{SP} \phantom{SP} {\tt else} \\
\phantom{SP} \phantom{SP} - set $m:=M$ and {\tt goto} $(v2)$\\
$(v8)$ {\tt if} $res>\varepsilon \, \& \, \lceil \theta m\rceil \leq M$ {\tt then}\\
\phantom{SP} \phantom{SP} - set $m:=\lceil \theta m\rceil$, $res_0:=res$ and {\tt goto} $(v2)$ 

\item[$(vi)$] select the indexes $J=\{i:\,u_i>0\}$ and set $w:=u(J)$ and $Z:=X(J)$

\item {\tt output}: the weights $\{w_k\}$ and nodes $\{Z_k\}\subset X$ of a compressed QMC formula 
on $\mathcal{J}$ with moment-matching residual $res$

\end{itemize}
\vskip0.3cm

Now, some observations on delicate aspects are in order. 
Step $(iii)$ is a key point in the case 
of surface integration. 
As for the starting polynomial basis, for conditioning problems we adopt the product Chebyshev total-degree basis of the smaller bounding box say $[a_1,b_1]\times[a_2,b_2]\times [a_3,b_3]\supset X$, namely 
$$
p_j(x,y,z)=T_{\alpha_1(j)}(\sigma_1(x))\cdot T_{\alpha_2(j)}(\sigma_2(y))\cdot T_{\alpha_3(j)}(\sigma_3(z)),\;\;j=1,\dots,\mathcal{V}\;,
$$
$$
\sigma_i:\,[a_i,b_i]\mapsto [-1,1]\;,\;\;\sigma_i(t)=\frac{2t-b_i-a_i}{b_i-a_i}\;,\;\;i=1,2,3\;,
$$
where $j\mapsto \alpha(j)$ corresponds to the graded lexicographical ordering of the 3-indexes $\alpha=(\alpha_1,\alpha_2,\alpha_3)$, $0\leq \alpha_1+\alpha_2+\alpha_3\leq n$.

Moreover, we recall that $dim(\mathbb{P}_n^3(X))$ is simply the rank of the corresponding rectangular Vandermonde-like matrix $C$. In step $(iii4)$, instead, we work with the principal square submatrix $C_{\mathcal{V}}$. As already observed in Section 2, with $\mathcal{V}\geq dim(\mathbb{P}_n^3(\mathcal{J}))$ uniformly distributed points on $\mathcal{J}$, the probability that such a rank be lower than 
$dim(\mathbb{P}_n^3(\mathcal{J}))$ is null, so that ``almost-surely'' Wilhelmsen theorem applies. In Matlab, one can use directly the built-in function {\tt rank} based on an economy-size version of SVD. 
Notice that we are using a numerical rank (obtained by discarding the singular values below a tolerance close to machine precision), not the true rank. Nevertheless, dealing with polynomials restricted to $X$ this is numerically equivalent to work, up to very small errors, with the true polynomial space.
We stress that when $\mathcal{V}\ll M$, using $C_{\mathcal{V}}$ instead of $C$ gives experimentally a substantial speed-up to the rank computation, by a factor roughly of the order of $M/\mathcal{V}$.

The polynomial basis selection, i.e. the determination of a set of linearly independent polynomials on $\mathcal{J}$ within the starting basis, is performed in $(iii4)$ by a QR factorization with column pivoting of the Chebyshev-Vandermonde matrix $C_{\mathcal{V}}$ (again, an economy-size version can be used in Matlab that produces only the first $N$ columns of $Q$ and a column permutation vector). In such a way we 
select a polynomial basis of $\mathbb{P}_n^3(X)$ by the first $N$ components $\pi_1,\dots,\pi_N$ of the column permutation, say $(f_1,\dots,f_N)=(p_{\pi_1},\dots,p_{\pi_N})$.

We can now turn to the second key step of the algorithm, that is the extraction of a compressed 
QMC formula in $(v)$. As already observed, this is based 
on Wilhelmsen theorem, using just $X_M$ as extraction set, in a ``bottom-up'' fashion. This procedure avoids working directly on the complete matrix $V_M$ (cf. (\ref{big-NNLS})), as done instead 
in other previous approaches to QMC compression 
like \cite{DME18}, cf. also the discussion in \cite{BDME16,ESV22}. Indeed, the overall number of points, i.e. of rows of $V_M$, can be huge, up to the order of $10^5-10^6$. 
In practice, we proceed along increasingly dense subsequences of the overall sequence, 
solving the corresponding NNLS problems and stopping when the relative moment-matching 
residual becomes sufficiently small. 

To this purpose the classical Lawson-Hanson iterative method turns out to be a good choice, since it automatically seeks a sparse solution with a number of nonzeros not exceeding $N$. The method is implemented in most numerical programming environments, e.g. in Matlab by the built-in function {\tt lsqnonneg}. On the other hand, there are improvements of the algorithm, cf. for example 
\cite{Sl} for a survey, and the recent implementation named LHDM based on the concept of ``Deviation Maximization'' instead of ``column pivoting'' for the underlying QR factorizations, 
cf. \cite{DODM22,DMV20}. Indeed, in the present framework we have adopted LHDM, since it gives experimentally a speed-up of at least 2 with respect to {\tt lsqnonneg}.

In order to cope ill-conditioning of the matrices used in the sequence of NNLS problem, that worsens increasing the degree, we perform an orthogonalization of $V_m$ by QR factorization, that corresponds to work with the discrete orthogonal basis $(f_1,\dots,f_N)R_m^{-1}$. Such a basis is orthogonal with respect the counting measure supported at $X$, 
i.e. with respect to the discrete scalar product $\langle f,g \rangle_{X_m}=\sum_{i=1}^m f(P_i) g(P_i)$. Consequently, 
the original QMC moments have to be modified as in $(v4)$. 

It should be stressed that, due to the inherited ill-conditioning of the triangular factor $R_m$ by $V_m$, that increases with the degree, explicit inversion of $R_m$ in $(v4)$ is avoided by solving  linear systems via Gaussian elimination with row pivoting (that is in Matlab simply by applying the backslash operator).

We also notice that  the complete matrix $V_M$ is used only to compute the QMC moments in $(iv)$, unless $(v7a)$ has to be followed due to a residual decrease factor below the required threshold. Such a phenomenon turns out to occur seldom with high degrees and strong ill-conditioning. In such a case, computation of $A_M=V_MR_m^{-1}$ becomes the computational bulk slowing down the whole process.

\section{Numerical tests and demos}
In order to show the effectiveness of the bottom-up compression procedure 
of QMC surface integration, we present some numerical tests, where  
we compare ``Caratheodory-Tchakaloff'' compression of multivariate discrete measures as implemented 
in the general-purpose package {\em dCATCH} \cite{DMV20-2}, with the bottom-up approach described above. The Matlab codes and demos, collected in a package named {\em Qsurf}, are freely available at 
\cite{ESV23-2}.

In all the tests we have set the 
parameters of the algorithm to $\varepsilon=10^{-10}$, $\theta=2$, $\tau=10$, and $m$ has been inizialized to $2N$.
The tests have been performed with a CPU AMD Ryzen 5 3600 with 48 GB of RAM, running Matlab R2022a.

\subsection{Sphere region}
\begin{figure}
\centering
\includegraphics[height=3in]{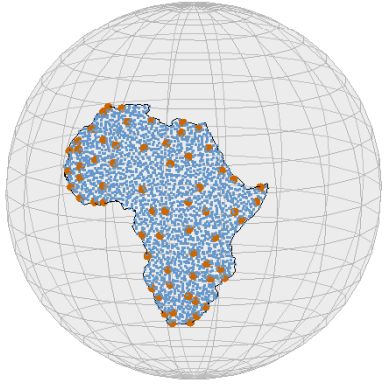}
\caption{64 compressed QMC points (red) at exactness degree $n=7$, extracted from about 2400 mapped Halton points (blue) on the surface of a spherical polygon approximating  continental Africa.}
\label{figure:africa}
\end{figure}

In the first example we consider a large region $\mathcal{J}$ of the sphere, namely a {\em spherical polygon} (a polygon whose vertices are on the sphere and whose sides are great circle arcs) representing an approximation of continental Africa (see Fig. \ref{figure:africa}). In this case it is convenient to choose a spherical cap (say $\mathcal{C}\supset \mathcal{J}$) centered at the polygon centroid as bounding surface, $\mathcal{S}=\mathcal{C}$, and we can apply a rotation to the sphere in such a way that the centroid is at the north pole (this does not clearly affect surface integration on the region). 

The indicator function of $\mathcal{J}$ can be easily implemented by stereographic projection from the south pole on the tangent plane at the north pole,  
that generates a planar polygon for which the Matlab {\tt inpolygon} works quite efficiently. Observe that this procedure can be applied to any rotated spherical polygon 
that does not contain the south pole. 

Then, we can parametrize the polar cap by the area-preserving 
map (i.e., $\|\partial_u\Psi\times \partial_v\Psi\|_2=1$)
\begin{equation} \label{area-pres}
\Psi(u,v)=r(\sqrt{1-u^2}\,\cos(v),\sqrt{1-u^2}\,\sin(v),u)\;,
\end{equation}
$(u,v)\in D=(c,1)\times (0,2\pi)$, 
where $r$ is the sphere radius and $c$ is the $z$-quote of the cap boundary (in practice, working with the open rectangle $D$ we loose the Greenwich $0$-meridian arc cutting the cap, that 
has null surface measure and thus surface integration is not affected). 

Now, starting from low-discrepancy points in $D$, e.g. Halton points, 
we get low-discrepancy points on the cap $\mathcal{S}$ and finally on the  
spherical polygon $\mathcal{J}$. On the other hand, Proposition 1 substantially applies since the map $\Psi$ is analytic and regular on $D$ (see also Remarks 1-2), and hence we can resort to  the bottom-up algorithm {\em Qsurf} in order to compress QMC integration on a huge number $M$ of mapped low-discrepancy points in $\mathcal{J}$.
To the purpose of illustration, in Figure {\ref{figure:africa}} we show the distribution of $64$ compressed QMC points extracted from about 2400 Halton points, still matching the QMC moments on  $\mathcal{J}$, up to degree 7. 

In Table {\ref{table:africa}} we report the results obtained by applying the QMC compression with more than one million points on the spherical polygon, taking degrees $n=3,6,9,12,15$, and accepting (\ref{res-tol}) with a tolerance $\epsilon=10^{-10}$. 
In particular, we display the cardinalities and compression ratios, the cpu-times for the construction of
the low-discrepancy sequence (cpu Halton seq.) and those for the computation of the compressed rules. 

The advantage of the new approach is two-fold, since in all the tests an inferior cputime with respect to {\em dCATCH} is required to determine the compressed rule and, differently from {\em{dCATCH}}, the solution of (\ref{NNLS}) always satisfies the moment residual criterion (\ref{res-tol}). In addition, less memory is necessary due to the inherent structure of the bottom-up approach, which works on much smaller matrices.

Finally, in Table \ref{table:africafunctions}, we approximate the integrals  $\int_{\mathcal{J}} g_k \,d\sigma$ on three test functions, namely setting  $P=(x,y,z)$
\begin{eqnarray}
g_1(P) &=& \exp(-\|P-P_0\|_2) \label{g1}\\
g_2(P) &=& \cos(x+y+z) \label{g2}\\
g_3(P) &=& \|P-P_0\|^5_2 \label{g3}
\end{eqnarray}
$P_0$ being the centroid of the spherical polygon $\mathcal{J}$.
The reference values of the integrals  hace been computed by a QMC rule with very high cardinality (more than 20 million points).  We display the relative errors of the QMC rule with more than one million points and of the two proposed compressions. As expected from estimate (\ref{equiv}), by increasing the QMC moment-matching degree the errors tend to stabilize around the underlying QMC error.

\begin{table}[ht]
\begin{center}
{\footnotesize
\begin{tabular}{| c | c | c | c | c | c |}
\hline
deg & 3 & 6 & 9 & 12 & 15 \\
\hline \hline
card. $QMC$ & \multicolumn{5}{c|}{$M=$ 1,184,341}  \\
\hline
card. $dCATCH$ & 16 & 49 & 98 & 165 & 239 \\
card. $\mbox{\em{Qsurf}}$    & 16 & 49 & 100 & 169 & 256 \\
compr. ratio & 7.4e+04 & 2.4e+04 & 1.1e+04 & 7.0e+03 & 4.6e+03 \\
\hline \hline  
cpu Halton seq.    & \multicolumn{5}{c|}{4.53e+01s} \\ 
\hline 
cpu $dCATCH$ & 4.1e+00s  & 1.5e+01s    & 4.6e+01s   & 1.3e+02s & 3.1e+02s  \\
\hline
cpu $\mbox{\em{Qsurf}}$ & 3.8e-01 & 1.2e+00s  & 3.1e+00s  & 6.4e+00s &1.3e+01   \\
speed-up & 10.8 & 12.5 & 14.8 & 20.3 & 23.8 \\
\hline \hline
mom. resid. $dCATCH$ & 4.3e-12 & 4.3e-12 & $\star$ 2.9e-04 & $\star$ 6.2e-04 & $\star$ 2.0e-03\\
\hline
mom. resid. $\mbox{\em{Qsurf}}$ & & & & & \\
iter. 1 & 3.7e-16 & 8.5e-01 & 2.9e+01 & 7.9e+01 & 1.4e+01\\
iter. 2 &          & 6.5e-02 & 1.7e-04 & 3.8e-02 & 7.6e-01\\
iter. 3 &          & 6.3e-16 & 1.1e-15 & 1.3e-15 & 1.7e-02\\
iter. 4 &          &          &          &          & 2.8e-15\\
\hline
\end{tabular}
}
\caption{\small{QMC compression by with more than one million points on a spherical polygon approximating  continental Africa.}}
\label{table:africa}
\end{center}
\end{table}

\begin{table}[h]
\begin{center}
{\footnotesize
\begin{tabular}{| c | c | c | c | c | c |}
\hline
deg & 3 & 6 & 9 & 12 & 15 \\
\hline \hline
$E^{QMC}(g_1)$  & \multicolumn{5}{c|}{ 3.0e-05}  \\ 
\hline 
$E^{dCATCH}(g_1)$ & 1.0e-03 & 3.1e-05 & 2.7e-05 & 3.2e-05 & 1.7e-05  \\
$E^{\mbox{\em{Qsurf}}}(g_1)$    & 1.2e-04 & 3.0e-05 & 3.0e-05 & 3.0e-05 & 3.0e-05 \\
 \hline \hline
$E^{QMC}(g_2)$  & \multicolumn{5}{c|}{ 1.5e-05}  \\ 
\hline  
$E^{dCATCH}(g_2)$ & 8.9e-05 & 1.5e-05 & 2.5e-05 & 4.9e-07 & 3.7e-06 \\
$E^{\mbox{\em{Qsurf}}}(g_2)$    & 1.4e-05 & 1.5e-05 & 1.5e-05 & 1.5e-05 & 1.5e-05 \\
 \hline \hline
$E^{QMC}(g_3)$  & \multicolumn{5}{c|}{8.6e-04 }  \\ 
\hline 
$E^{dCATCH}(g_3)$ & 2.4e-02 & 1.3e-03 & 7.8e-04 & 8.2e-04 & 5.8e-04 \\
$E^{\mbox{\em{Qsurf}}}(g_3)$    & 2.3e-02 & 7.7e-04 & 8.3e-04 & 8.6e-04 & 8.6e-04 \\
\hline 
\end{tabular}
}
\caption{\small{Relative integration errors 
for the three test functions (\ref{g1})-(\ref{g3}) on a spherical polygon approximating continental Africa, by means of QMC, {\em{dCATCH}} and {\em{Qsurf}}.
}}
\label{table:africafunctions}
\end{center}
\end{table}

\subsection{Torus region}
\begin{figure}
    \centering
    \subfigure[]
        {\includegraphics[height=1.81in]{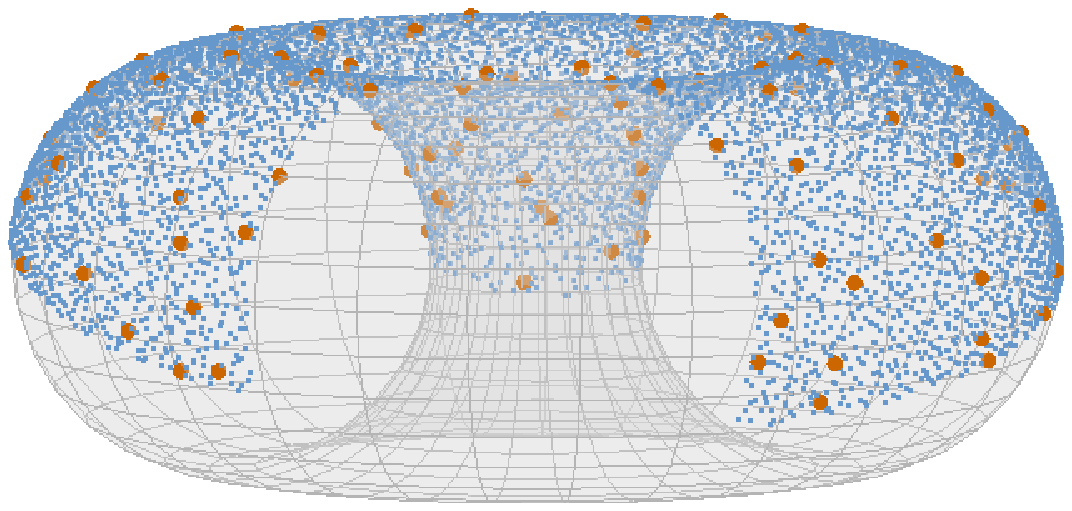}}
     ~ 
    \subfigure[]
        {\includegraphics[height=1.82in]{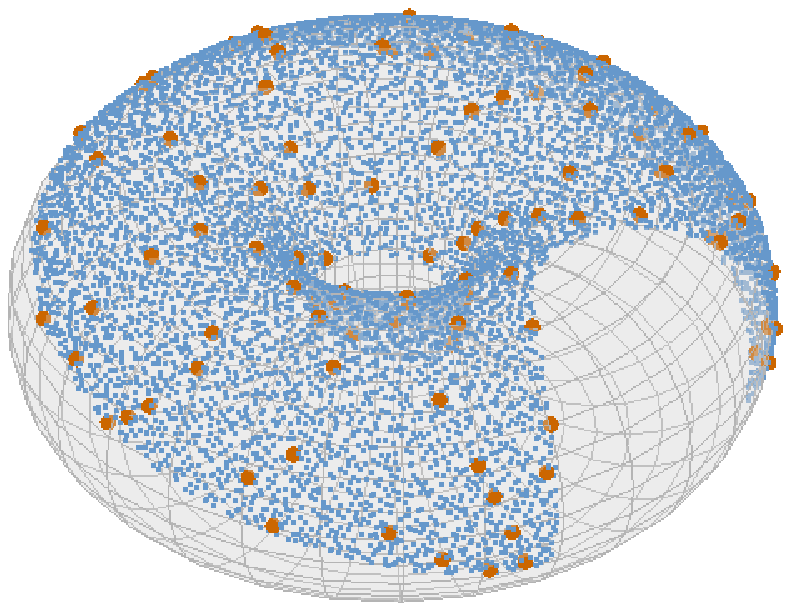}}
    \caption{100 compressed QMC points (red) at exactness degree $n=7$, extracted from about 8000 mapped Halton points (blue) by rejection sampling on a torus region, determined by a cutting ball and plane (view from different perspectives).}
    \label{figure:torus}
\end{figure}

The second example concerns surface integration on a region $\mathcal{J}$ of a torus $\cal{T}$, corresponding to a section by a plane, excluding the points that are internal 
to a ball intersecting the torus; see Fig. \ref{figure:torus}. 
In particular, we consider the torus with center $(0,0,0)$ and radii $r=2$, $R=3$, cut by the ball $B((0,4,0),\sqrt{6})$ and the upper half-space of ${\mathbb{R}}^3$ w.r.t. the plane of equation $-x/4+y+4z=0$. 

In this case it is not straightforward to apply a standard integrator, since one should track the domain 
$\Psi^{-1}(\mathcal{J})$ in $D$ and then apply there a suitable cubature rule. On the contrary, QMC integration can be more easily constructed by rejection sampling in standard toroidal coordinates (here the bounding surface $\mathcal{S}=\cal{T}$ is the whole torus) 
\begin{equation} \label{tor-coord}
\Psi(u,v)=((R+r\cos(u))\cos(v),(R+r\cos(u))\sin(v),r\sin(u))\;,
\end{equation}
$(u,v)\in D=(0,2\pi)\times (0,2\pi)$, where $R$ and $r$ are the big and small torus radii respectively, and $\|\partial_u\Psi\times \partial_v\Psi\|_2=r(R+r\cos(v))$. Observe that considering the open rectangle $\Omega$ we loose the possible intersection of $\mathcal{J}$ with two circles, that have null surface measure and do not affect surface integration.  
Moreover, the indicator function of $\mathcal{J}$ can be implemented by the simple inequalities that describe an half-space determined by the cutting plane, and the interior of the ball. Again, the map $\Psi$ is analytic and regular so that Proposition 1 with Remarks 1-2 applies and algorithm {\em Qsurf} can be used. In Figure {\ref{figure:torus}} we show the distribution of $64$ compressed QMC points, extracted from about 8000 mapped Halton points after selection by rejection sampling w.r.t. the surface measure density, still matching the QMC moments on $\mathcal{J}$ up to degree 7.

In Table {\ref{table:torus}} we again report the results obtained by applying QMC compression with more than one million points on the region $\mathcal{J}$. 
As for the spherical polygon, we consider degrees $n=3,6,9,12,15$, accepting (\ref{res-tol}) with a tolerance $\epsilon=10^{-10}$.
In all the tests an inferior cputime is required by {\em{Qsurf}} to determine the compressed rule and, while {\em{dCATCH}} fails for degree $n=15$. Moreover, the solution by the new approach to (\ref{NNLS}) always satisfies the moment residual criterion (\ref{res-tol}). 

Lastly, in Table \ref{table:torusfunctions} we approximate the value of $\int_{\mathcal{S}} g_k\,d\sigma$, $k=1,2,3$, with the same functions defined in ({\ref{g1}})-(\ref{g3}) and $P_0=(0,-3,2)$. The reference values of the integrals have been computed by means of a QMC rule with very high cardinality (more than 20 million points). We display the relative errors of the QMC rule with about one million points and of the two proposed compressions. 
Notice again that, as expected from estimate (\ref{equiv}), by increasing the QMC moment-matching degree the errors tend to stabilize around the underlying QMC error.

\begin{table}[ht]
\begin{center}
{\footnotesize
\begin{tabular}{| c | c | c | c | c | c |}
\hline
deg & 3 & 6 & 9 & 12 & 15 \\
\hline \hline
card. $QMC$ & \multicolumn{5}{c|}{$M=$ 1,006,200}  \\
\hline
card. $dCATCH$ & 20 & 74 & 164 & 290 & 450 \\
card. ${\mbox{\em Qsurf}}$    & 20 & 74 & 164 & 290 & 452 \\
compr. ratio & 5.0e+04 & 1.3e+04 & 6.1e+03 & 3.5e+03 & 2.2e+03 \\
\hline \hline  
cpu Halton seq. & \multicolumn{5}{c|}{1.0e+01s} \\ 
\hline 
cpu $dCATCH$ & 2.8e+00s  & 1.6e+01s   & 4.4e+01s   & 1.2e+02s  & 3.0e+02s   \\
\hline
cpu $\mbox{\em{Qsurf}}$ & 2.7e-01s  & 9.9e-01s  & 2.9e+00s  & 6.3e+00s &2.2e+01s   \\
speed-up & 10.4 & 16.2 & 15.2 & 19.0 & 13.6 \\
\hline \hline
mom. resid. $dCATCH$ & 1.2e-11 & 1.2e-11 & 1.2e-11 & 1.2e-11 & $\star$ 9.1e-07\\
\hline
mom. resid. $\mbox{\em{Qsurf}}$ & & & & & \\
iter. 1 & 3.0e-16 & 8.9e-01 & 1.3e+00 & 6.4e+00 & 2.5e+01\\
iter. 2 &          & 1.1e-15 & 1.9e-15 & 2.6e-01 & 1.3e-01\\
iter. 3 &          &          &          & 3.3e-15 & 4.5e-15\\
\hline
\end{tabular}
}
\caption{\small{QMC compression with more than one million points on the torus region in Fig. \ref{figure:torus}.}}
\label{table:torus}
\end{center}
\end{table}

\begin{table}[h]
\begin{center}
{\footnotesize
\begin{tabular}{| c | c | c | c | c | c |}
\hline
deg & 3 & 6 & 9 & 12 & 15 \\
\hline \hline
$E^{QMC}(g_1)$  & \multicolumn{5}{c|}{ 1.7e-04}  \\ 
\hline 
$E^{dCATCH}(g_1)$ & 3.5e-01 & 1.2e-02 & 2.5e-03 & 2.2e-04 & 2.2e-04 \\
$E^{\mbox{\em{Qsurf}}}(f_1)$    & 5.5e-01 & 6.5e-02 & 2.4e-03 & 5.4e-04 & 1.5e-04 \\
 \hline \hline
$E^{QMC}(g_2)$  & \multicolumn{5}{c|}{ 2.4e-04}  \\ 
\hline  
$E^{dCATCH}(g_2)$ & 3.5e-01 & 2.5e-01 & 7.2e-03 & 1.3e-04 & 2.4e-04 \\
$E^{\mbox{\em{Qsurf}}}(f_2)$    & 1.7e+00 & 1.3e-01 & 1.5e-03 & 1.8e-04 & 2.4e-04 \\
 \hline \hline
$E^{QMC}(g_3)$  & \multicolumn{5}{c|}{5.2e-06 }  \\ 
\hline 
$E^{dCATCH}(f_3)$ & 4.3e-03 & 2.3e-06 & 5.2e-06 & 5.2e-06 & 5.2e-06 \\
$E^{\mbox{\em{Qsurf}}}(g_3)$    & 8.0e-03 & 3.1e-06 & 5.2e-06 & 5.2e-06 & 5.2e-06 \\
\hline 
\end{tabular}
}
\caption{\small{Relative errors for the three test functions (\ref{g1})-(\ref{g3}) on the torus region of Fig. {\ref{figure:torus}},  by means of QMC, {\em{dCATCH}} and {\em{Qsurf}}.
}}
\label{table:torusfunctions}
\end{center}
\end{table}

\subsection{Cartesian graph}
In the third example we consider as a regular surface $\mathcal{S}$ the Cartesian graph of an analytic function, namely the popular Franke's surface, which is the graph of a linear combination of Gaussians 
$$
F(u,v)=\frac{3}{4}\,e^{-\frac{1}{4}\,((9u-2)^2+(9v-2)^2)}+\frac{3}{4}\,e^{-\frac{1}{49}\,((9u+1)^2+(9v+1)^2)}
$$
\begin{equation} \label{franke}
+\frac{1}{2}\,e^{-\frac{1}{4}\,((9u-7)^2+(9v-3)^2)}-\frac{1}{5}\,e^{-((9u-4)^2+(9v-7)^2)}\;,
\end{equation}
$(u,v)\in D=(0,1)\times (0,1)$.

\begin{figure}
    \centering
    \subfigure[]
        {\includegraphics[height=1.81in]{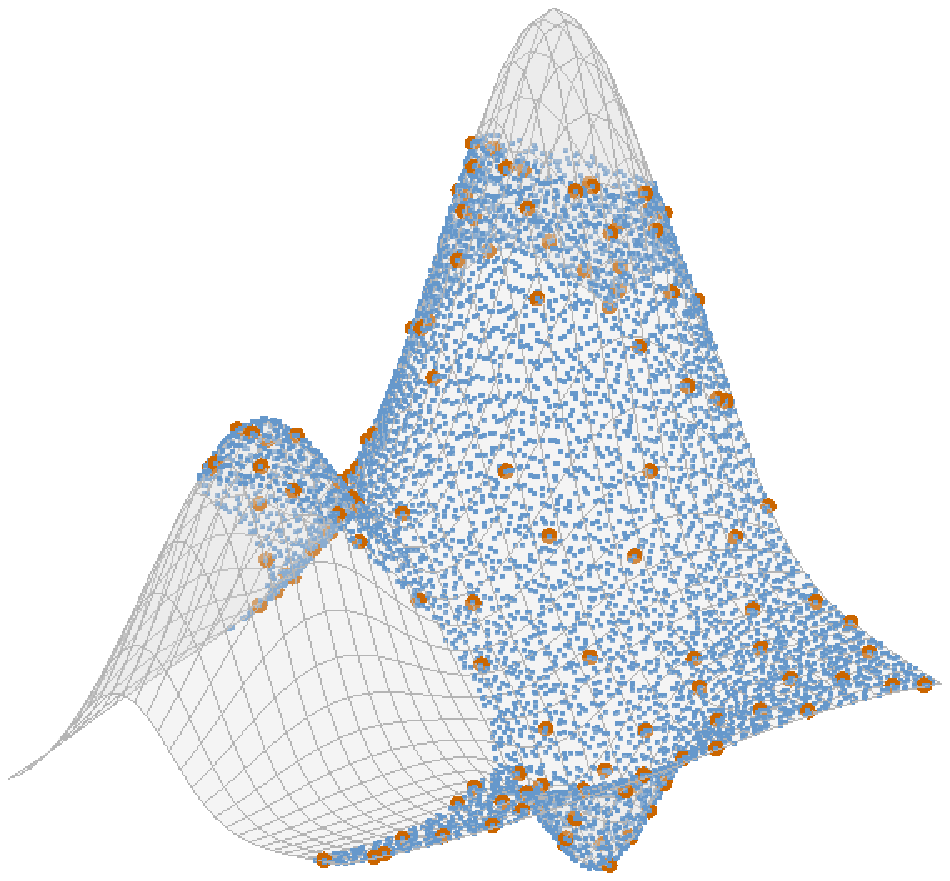}}
     ~ 
    \subfigure[]
        {\includegraphics[height=1.82in]{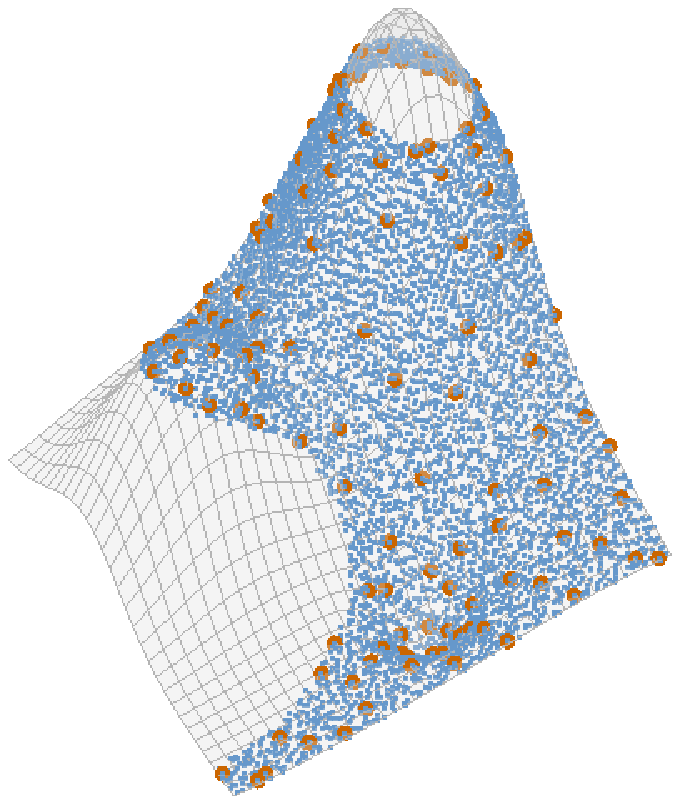}}
    \caption{120 compressed QMC points (red) at exactness degree $n=7$, extracted from about 6500 mapped Halton points (blue) by rejection sampling on a Franke's surface region, determined by a cutting ball and plane (view from different perspectives).}
    \label{franke1}
\end{figure}

\begin{figure}
    \centering
    \subfigure[]
        {\includegraphics[height=1.81in]{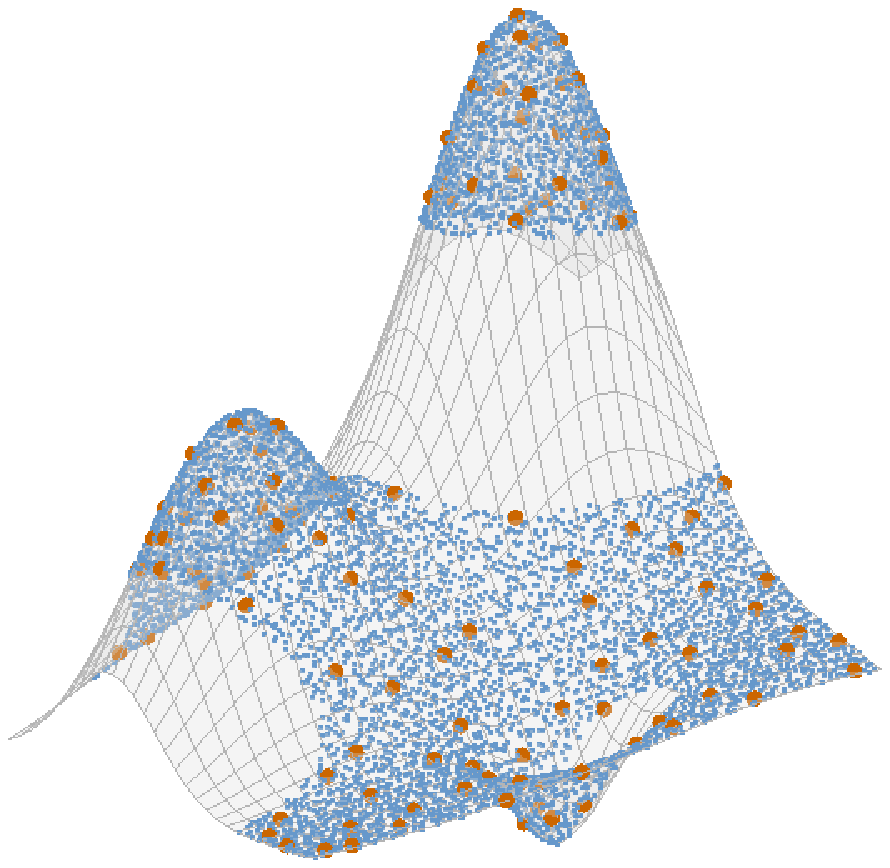}}
     ~ 
    \subfigure[]
        {\includegraphics[height=1.82in]{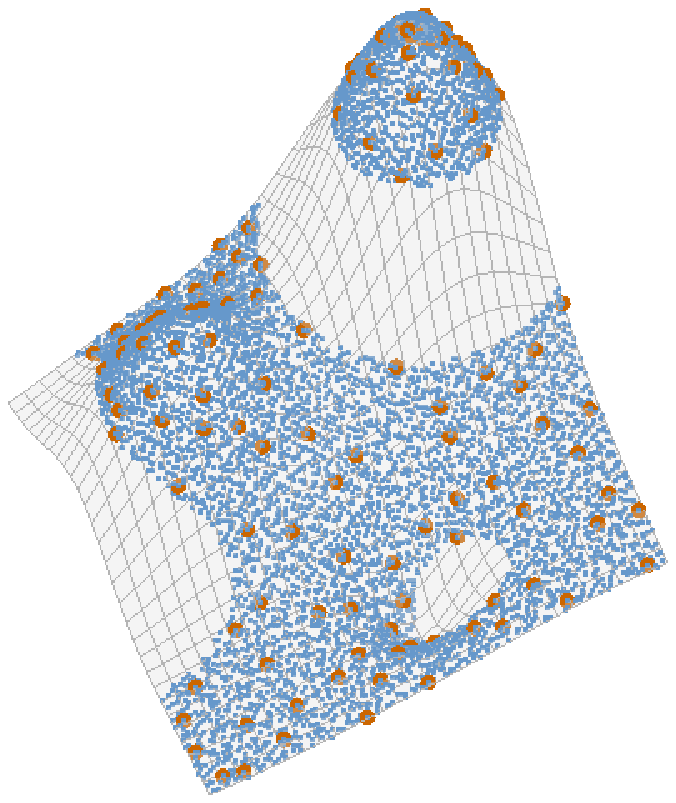}}
    \caption{120 compressed QMC points (red) at exactness degree $n=7$, extracted from about 6500 mapped Halton points (blue) by rejection sampling on a Franke's surface disconnected region, determined by three cutting balls (view from different perspectives).}
    \label{franke2}
\end{figure}

We take two regions of such a surface, the first determined by a cutting ball and plane, whereas the second is a disconnected one determined by three cutting balls; see Figs. \ref{franke1} and \ref{franke2}.  
Again, the map $\Psi$ is analytic and regular, since $\|\partial_u\Psi\times \partial_v\Psi\|_2=\sqrt{1+(\partial_u F)^2+(\partial_v F)^2}$, so that Proposition   1 with Remarks 1-2 applies and algorithm {\em Qsurf} can be used.

The numerical tests are collected in Tables 5-8, and show results that are in line with those of the previous examples, apart from the fact that the numerically determined dimension of the trivariate polynomial spaces does not collapse on the surface (at least up to degree 9). This is expected since  Franke's surface is a transcendental, i.e. not algebraic, surface. Notice in particular that 
at degrees 9, 12, 15, {\em dCATCH} fails to reach the 
required residual tolerance, whereas {\em Qsurf} always 
succeeds in at most 4-5 iterations.

\begin{table}[ht]
\begin{center}
{\footnotesize
\begin{tabular}{| c | c | c | c | c | c |}
\hline
deg & 3 & 6 & 9 & 12 & 15 \\
\hline \hline
card. $QMC$ & \multicolumn{5}{c|}{$M=$ 1,293,600}  \\
\hline
card. $dCATCH$ & 20 & 84 & 212 & 407 & 586 \\
card. $\mbox{\em{Qsurf}}$    & 20 & 84 & 220 & 442 & 701 \\
compr. ratio & 6.5e+04 & 1.5e+04 & 5.9e+03 & 2.9e+03 & 1.8e+03 \\
\hline \hline  
cpu Halton seq. & \multicolumn{5}{c|}{1.8e+00s} \\ 
\hline 
cpu $dCATCH$ & 3.6e+00s  & 2.1e+01s  & 5.7e+01s   & 1.8e+02s & 4.9e+02s   \\
\hline
cpu $\mbox{\em{Qsurf}}$ &  3.5e-01s   &  1.4e+00s   &  3.8e+00s   &  2.0e+01s  &  2.7e+01s   \\
speed-up & 10.3  &  15.0 & 15.0 & 9.0 & 18.1 \\
\hline \hline
mom. resid. $dCATCH$ & 7.6e-12 & 7.6e-12 & $\star$ 8.9e-04 & $\star$ 2.9e-03 & $\star$ 5.9e-03\\
\hline
mom. resid. $\mbox{\em{Qsurf}}$ & & & & & \\
iter. 1 & 1.9e-16 & 5.5e-01 & 1.5e+00 & 1.4e+01 & 3.4e+01 \\
iter. 2 &          & 1.0e-15 & 1.8e-15 & 3.9e-01 & 1.9e+00 \\
iter. 3 &          &          &          & 1.2e-02 & 4.2e-15 \\
iter. 4 &          &          &          & 2.5e-15 &          \\
\hline
\end{tabular}
}
\caption{\small{QMC compression with more than one million points on 
on the Franke's surface region in Fig. \ref{franke1}.}}
\label{3balls}
\end{center}
\end{table}

\begin{table}[ht]
\begin{center}
{\footnotesize
\begin{tabular}{| c | c | c | c | c | c |}
\hline
deg & 3 & 6 & 9 & 12 & 15 \\
\hline \hline
card. $QMC$ & \multicolumn{5}{c|}{$M=$ 1,305,444}  \\
\hline
card. $dCATCH$ & 20 & 84 & 212 & 405 & 612 \\
card. $\mbox{\em{Qsurf}}$    & 20 & 84 & 220 & 448 & 735 \\
compr. ratio & 6.5e+04 & 1.6e+04 & 5.9e+03 & 2.9e+03 & 1.8e+03 \\
\hline \hline  
cpu Halton seq. & \multicolumn{5}{c|}{1.74e+00s} \\ 
\hline 
cpu $dCATCH$ & 3.8e+00 s  & 2.2e+01   & 5.2e+01  & 1.9e+02s & 5.0e+02s   \\
\hline
cpu $\mbox{\em{Qsurf}}$ &  3.6e-01    &  1.4e+00   &  9.0e+00    &  3.9e+01 &  6.3e+01    \\
speed-up & 10.6  & 15.7  & 5.8  & 4.9  & 7.9 \\
\hline \hline
mom. resid. $dCATCH$ & 1.8e-12 & 1.8e-12 & $\star$ 8.2e-04 & $\star$ 2.4e-03 & $\star$ 3.4e-03\\
\hline
mom. resid. $\mbox{\em{Qsurf}}$ & & & & & \\
iter. 1 & 6.0e-16 & 7.4e-02 & 7.7e-01 & 3.2e+01 & 3.4e+01 \\
iter. 2 &          & 1.2e-15 & 2.4e-01 & 4.0e-01 & 6.1e+00 \\
iter. 3 &          &          & 2.4e-01 & 8.5e-02 & 6.1e+00 \\
iter. 4 &          &          & 1.6e-11 & 8.5e-02 & 9.2e-12 \\
iter. 5 &          &          &          & 7.9e-12 &          \\
\hline
\end{tabular}
}
\caption{\small{QMC compression with more than one million points on the on the Franke's surface disconnected region in Fig. \ref{franke2}.}}
\label{3balls}
\end{center}
\end{table}

\begin{table}[h]
\begin{center}
{\footnotesize
\begin{tabular}{| c | c | c | c | c | c |}
\hline
deg & 3 & 6 & 9 & 12 & 15 \\
\hline \hline
$E^{QMC}(g_1)$  & \multicolumn{5}{c|}{ 1.2e-05}  \\ 
\hline 
$E^{dCATCH}(g_1)$ & 4.4e-03 & 1.1e-05 & 1.5e-06 & 8.1e-06 & 4.6e-05 \\
$E^{\mbox{\em{Qsurf}}}(g_1)$    & 6.7e-04 & 9.6e-06 & 1.2e-05 & 1.2e-05 & 1.2e-05 \\
 \hline \hline
$E^{QMC}(g_2)$  & \multicolumn{5}{c|}{ 3.0e-07}  \\ 
\hline  
$E^{dCATCH}(g_2)$ & 7.8e-05 & 2.9e-07 & 3.6e-05 & 1.2e-04 & 6.3e-05 \\
$E^{\mbox{\em{Qsurf}}}(g_2)$    & 4.0e-05 & 3.0e-07 & 3.0e-07 & 3.0e-07 & 3.0e-07 \\
 \hline \hline
$E^{QMC}(g_3)$  & \multicolumn{5}{c|}{ 6.0e-05}  \\ 
\hline 
$E^{dCATCH}(g_3)$ & 1.4e-01 & 3.8e-05 & 1.1e-04 & 2.8e-05 & 1.2e-04 \\
$E^{\mbox{\em{Qsurf}}}(g_3)$    & 2.2e-02 & 1.1e-04 & 6.1e-05 & 6.0e-05 & 6.0e-05 \\
\hline 
\end{tabular}
}
\caption{\small{Relative integration errors for the three test functions (\ref{g1})-(\ref{g3}) on the Franke's surface region of Fig. \ref{franke1}, by means of QMC, $dCATCH$, {\em Qsurf}.}}
\label{3balls-errors}
\end{center}
\end{table}

\begin{table}[h]
\begin{center}
{\footnotesize
\begin{tabular}{| c | c | c | c | c | c |}
\hline
deg & 3 & 6 & 9 & 12 & 15 \\
\hline \hline
$E^{QMC}(g_1)$  & \multicolumn{5}{c|}{1.3e-06 }  \\ 
\hline 
$E^{dCATCH}(g_1)$ & 4.6e-03 & 2.3e-06 & 1.7e-05 & 8.1e-06 & 2.1e-05 \\
$E^{\mbox{\em{Qsurf}}}(g_1)$    & 2.1e-03 & 9.7e-07 & 1.3e-06 & 1.3e-06 & 1.3e-06 \\
 \hline \hline
$E^{QMC}(g_2)$  & \multicolumn{5}{c|}{ 6.6e-05}  \\ 
\hline  
$E^{dCATCH}(g_2)$ & 5.3e-05 & 6.6e-05 & 6.7e-05 & 9.1e-05 & 4.8e-05 \\
$E^{\mbox{\em{Qsurf}}}(g_2)$    & 1.1e-04 & 6.6e-05 & 6.6e-05 & 6.6e-05 & 6.6e-05 \\
 \hline \hline
$E^{QMC}(g_3)$  & \multicolumn{5}{c|}{ 1.7e-06}  \\ 
\hline 
$E^{dCATCH}(g_3)$ & 2.3e-01 & 4.1e-05 & 1.6e-04 & 6.3e-05 & 8.1e-05 \\
$E^{\mbox{\em{Qsurf}}}(g_3)$    & 1.1e-01 & 4.0e-05 & 1.3e-06 & 1.7e-06 & 1.7e-06 \\
\hline 
\end{tabular}
}
\caption{\small{Relative integration errors for the three test functions (\ref{g1})-(\ref{g3}) on the Franke's surface disconnected region of Fig. \ref{franke2}, by means of QMC, $dCATCH$, {\em Qsurf}.}}
\label{3balls-errors}
\end{center}
\end{table}


\section{Software}

We have implemented and tested in Matlab all the described routines.

The demos {\tt{demo\_CQMC\_sphpoly}}, {\tt{demo\_CQMC\_torus}}, {\tt{demo\_CQMC\_franke}} illustrate the numerical experiments performed in the previous section. Their structure is essentially similar and can be modified to treat other subsets and/or parametric surfaces, adapting the function {\tt{pts\_domain}} to the new instance.
This corresponds to items ($i$) and ($ii$) of {\em{Algorithm Qsurf}}.

The routine {\tt{cqmc\_v2}} implements its remaining items from ($iii$) to ($vi$). To this purpose, the basis selection in ($iii$) is obtained by means of the function {\tt{dCHEBVAND\_v2}}, while the computation of a sparse solution in ($v5$) is achieved by an user's choice implementation of the Lawson-Hanson algorithm (namely, the Matlab built-in function {\tt{lsqnonneg}} or the alternative open-source codes {\tt{lawsonhanson}} and {\tt{LHDM}} proposed respectively in \cite{Sl} and \cite{DMV20}). Moreover, having in mind to compare algorithm Qsurf with previous approaches, we also provide the routine {\tt{dCATCH}} from \cite{DMSV20}, which implements Caratheodory-like compression via NNLS.

The open source software is available at \cite{ESV23-2}.

\section*{Acknowledgements}
Work partially
supported by the
DOR funds and the biennial project BIRD 192932
of the University of Padova, and by the INdAM-GNCS 
2022 Project ``Methods and software for multivariate integral models''.
This research has been accomplished within the RITA ``Research ITalian network on Approximation", the UMI Group TAA ``Approximation Theory and Applications" (G. Elefante, A. Sommariva) and the SIMAI Activity Group ANA\&A (A. Sommariva, M. Vianello).


\begin{thebibliography}{99}
 
\bibitem{BDME16} L. Bittante, S. De Marchi, G. Elefante, A new quasi-Monte Carlo technique based on nonnegative least-squares and approximate Fekete
points, Numer. Math. Theory Methods Appl. 9 (2016), 640--663.

\bibitem{BCCGST14} L. Brandolini, C. Choirat, L. Colzani, G. Gigante, R. Seri, L. Travaglini, 
Quadrature rules and distribution of points on manifolds, Ann. Sc. Norm. Super. Pisa Cl. Sci. (5)
Vol. XIII (2014), 889--923.

\bibitem{BSSW14} J.S. Brauchart, E,B. Saff, I.H. Sloan, R.S. Womersley, QMC designs: optimal order Quasi Monte Carlo Integration schemes on the sphere, Math. Comp. 83 (2014), 2821--2851. 

\bibitem{CLOS15} D.A. Cox, J. Little, D. O’Shea, Ideals, varieties, and algorithms, 4th edition, Springer, 2015.

\bibitem{D67} P.J. Davis,   
A construction of nonnegative approximate quadratures, 
Math. Comp. 21 (1967), 578--582.

\bibitem{DASV23} F. Dell'Accio, A. Sommariva, M. Vianello, Random sampling and unisolvent interpolation
by almost everywhere analytic functions,  arXiv:2303.14074. 

 \bibitem{DODM22} M. Dell'Orto, M. Dessole, F. Marcuzzi, The Lawson-Hanson Algorithm with Deviation Maximization: Finite Convergence and Sparse Recovery, Numer. Linear Algebra Appl., published online 13 January 2023. 

\bibitem{DME18} S. De Marchi, G. Elefante, Quasi-Monte Carlo integration on manifolds with mapped low-discrepancy points and greedy minimal Riesz
s-energy points, Appl. Numer. Math. 127 (2018), 110–-124.

\bibitem{DM21} M. Dessole, F. Marcuzzi, Deviation maximization for rank-revealing QR factorizations, Numer. Algorithms 91 (2022), 1047-1079.

\bibitem{DMV20}  M. Dessole, F. Marcuzzi, M. Vianello, Accelerating the Lawson-Hanson NNLS solver for large-scale Tchakaloff regression designs, Dolomites Res. Notes Approx. DRNA 13 (2020), 20--29.

\bibitem{DMV20-2} M. Dessole, F. Marcuzzi, M. Vianello, dCATCH: a numerical package for d-variate near G-optimal Tchakaloff regression via fast NNLS, 
MDPI-Mathematics 8(7) (2020) - Special Issue "Numerical Methods".

\bibitem{DMSV20} M. Dessole, F. Marcuzzi, M. Vianello, dCATCH: dCATCH: numerical package for d-variate discrete measure compression, near-optimal design and polynomial fitting - v1.1\\ 
\url{https://www.math.unipd.it/~marcov/dCATCH.html}.

\bibitem{DP10} J. Dick and F. Pillichshammer, Digital Nets and Sequences. Discrepancy Theory and Quasi-Monte Carlo Integration, Cambridge University
Press, Cambridge, 2010.

\bibitem{ESV22} G. Elefante, A. Sommariva, M. Vianello, CQMC: an improved code for low-dimensional Compressed Quasi-MonteCarlo cubature, 
Dolomites Res. Notes Approx. DRNA 15 (2022). 

\bibitem{ESV23} G. Elefante, A. Sommariva, M. Vianello, Compressed QMC volume and surface integration on union of balls, arXiv:2303.01460.

\bibitem{ESV23-2} G. Elefante, A. Sommariva, M. Vianello, Qsurf: a software package for compressed QMC integration on parametric surfaces (in Matlab)\\ 
\url{https://www.math.unipd.it/~alvise/software.html}.

\bibitem{H21} S. Hayakawa, Monte Carlo cubature construction, Jpn. J. Ind. Appl. Math. 38 (2021), 561-577.

\bibitem{LH95} 
 C.L. Lawson, R.J. Hanson, Solving least squares problems. Classics in Applied Mathematics 15, SIAM, Philadelphia, 1995.

\bibitem{L21} G. Legrain,  Non-Negative Moment Fitting Quadrature Rules for Fictitious Domain Methods, 
Comput. Math. Appl. 99 (2021), 270--291. 

\bibitem{LL12} C. Litterer, T. Lyons, High order recombination and an application to cubature on Wiener space, Ann. Appl. Probab. 22 (2012), 1301–1327. 

\bibitem{NO16} N. Nguyen, G. \"{O}kten, The acceptance-rejection method for low-discrepancy sequences, Monte Carlo Methods Appl. 22 (2016), 133--148.
 
\bibitem{PSV17}  
F. Piazzon, A. Sommariva, M. Vianello\/, Caratheodory-Tchakaloff Subsampling, Dolomites Res. Notes Approx. DRNA 10 (2017), 5--14.

\bibitem{P09} W. Pl\'{e}sniak, Multivariate Jackson Inequality, J. Comput. Appl. Math. 233 (2009), 815–820. 

\bibitem{R71} D.L. Ragozin, Constructive Polynomial Approximation on Spheres and Projective Spaces, Trans. Amer. Math. Soc. 162 (1971), 157--170. 
 
\bibitem{Sl} 
M. Slawski, Non-negative least squares: comparison of algorithms\\ \url{https://sites.google.com/site/slawskimartin}. 

\bibitem{SV15} 
 A. Sommariva, M. Vianello, Compression of multivariate discrete measures and applications,  Numer. Funct. Anal. Optim. 36 (2015), 1198--1223.

\bibitem{SV21} A. Sommariva, M. Vianello, Computing Tchakaloff-like cubature rules on spline curvilinear polygons, 
Dolomites Res. Notes Approx. DRNA 14 (2021), 1--11.

\bibitem{SV23} A. Sommariva, M. Vianello, Low-cardinality Positive Interior cubature on NURBS-shaped domains, 
BIT Numer. Math., published online 17 March 2023.

\bibitem{T57} V. Tchakaloff, Formules de cubatures m\'ecaniques \`a coefficients non n\'egatifs, (French), Bull. Sci. Math. 81 (1957), 123--134.
 
 \bibitem{Tche15} 
 M. Tchernychova, Caratheodory cubature measures. Ph.D. dissertation in Mathematics (supervisor: T. Lyons), University of Oxford, 2015. 
 
\bibitem{W76} 
D.R. Wilhelmsen, A Nearest Point Algorithm for Convex Polyhedral Cones and Applications to Positive Linear approximation, Math. Comp. 30 (1976), 48--57.

\bibitem{ZD14} H. Zhu, J. Dick, Discrepancy bounds for deterministic
acceptance-rejection samplers, Electron. J. Stat. 8 (2014), 678--707.

\end{thebibliography}
\end{document}